\documentclass
[12pt,a4paper]{article}
\usepackage[cp1251]{inputenc}
\usepackage{amssymb, amsfonts, amsmath, latexsym}
\usepackage[russian]{babel}
\begin{document}

\Large
\centerline{\bf A note-question on partitions of semigroups}\vspace{6 mm}

\normalsize\centerline{\bf  Igor Protasov, Ksenia Protasova}\vspace{6 mm}

{\bf Abstract.} Given a semigroup $S$  and an $n$-partition $\mathcal{P}$ of $S$, $n\in  \mathbb{N}$, do there exist $A\in \mathcal{P}$ and a subset $F$  of $S$ such that $S=F ^{-1} \{x \in S:  x  A \bigcap A\neq\emptyset\}$   and $|F |\leq n$?

We give an affirmative answer provided that either $S$ is finite or $n=2$.

\vspace{6 mm}

\vspace{3 mm}

\vspace{6 mm}

\centerline{\bf 1.	Introduction}
\vspace{3 mm}

In 1995, the first author asked the following question [3, Problem 13.44].

{\it Given a group $G$ and an   $n$-partition $\mathcal{P}$,  $n\in  \mathbb{N}$  of $G$, do there exist $A\in \mathcal{P}$ and a subset $F$ of $G$  such that $G=FAA ^{-1}$  and $|F |\leq n$?}

For the current state of this open problem see the survey [1].   We mention only that an answer  is positive if either $G$ is amenable (in particular, finite), or $n\leq 3$, or $x^{-1} A x =A$  for any $A\in  \mathcal{P}$  and $x\in  G$.  If $G$  is an arbitrary group and $\mathcal{P}$ is an $n$-partition of $G$ then one can choose $A$, $B\in \mathcal{P}$  and subsets $F$,  $H$  of $G$  such that $G=FAA^{-1}$, $|F|\leq  n!$ and $G=HBB^{-1} B$  and $|H|\leq n$.

In this note, we formulate a semigroup version of above question and give positive answer provided that either $S$ is finite or $n=2$.

For systematic exposition of Ramsey theory of semigroups see [2].

For a semigroup $S$, $a \in S$, $A \subseteq S$  and $B\subseteq S$, we use the standard notations
$$a^{-1} B= \{ x\in  S: ax\in  B\},  \ \  A^{-1} B=\bigcup  _{a\in A }  a^{-1} B.$$

We set $\Delta(A)=\{  x\in S: x A \bigcap A\neq\emptyset\}$  and observe that $\Delta (A)= \{x\in S: x^{-1} A \bigcap A\neq\emptyset\}$ and if $S$  is a group then $\Delta(A)=AA^{-1}.$

We suppose that $S^{-1} A=S$  and define a {\it covering number }
$$cov  A = min \{ |X|: X\subseteq  S, S = X^{-1} A\}.$$
If $S^{-1}A \neq  S$  then $cov A$ is not defined. Clearly, $cov A$   is defined if and only if $Sx \bigcap A \neq\emptyset $ for every $x\in S$.

Now we are ready for promised question.

{\it Given a semigroup $S$  and an $n$-partition $\mathcal{P}$  of $S$, does there exist $A \in \mathcal{P}$  such that $cov \Delta(A)\leq  n$?}
\vspace{3 mm}

\centerline{\bf 2.	Results}
\vspace{3 mm}

{\bf Theorem 1.} {\it For a semigroup $S$  and an $n$-partition $\mathcal{P}$  of $S$,  there exists $A\in \mathcal{P}$  such that
$cov \Delta (A)\leq 2 ^{2 ^{n-1}  -1}$. }

If $n=2$  then $cov \Delta(A) \leq 2$. In a personal communication, G. Bergman answered the question positively for $n=3$, and noticed that, we may suppose that a semigroup  $S$  is a monoid.
\vspace{3 mm}

{\bf Theorem 2.} {\it  For a finite semigroup $S$  and an $n$-partition $\mathcal{P}$  of $S$,  there exists $A \in \mathcal{P}$  such that $cov \Delta (A)\leq  n$ .}

\vspace{3 mm}

{\bf Theorem 3.} {\it If a subset $A$  of a semigroup $S$  contains either left or right zero then $cov \Delta (A)=1$.}
\vspace{3 mm}

\centerline{\bf 3.	Proofs}
\vspace{3 mm}

{\it Proof of Theorem 1.} We adopt arguments from [4, pp.  120-121] proving this theorem for groups.

We define a function $f: \mathbb{N}\times  \mathbb{N}\longrightarrow  \mathbb{N}$ by
$$f(1,m)=m  \ \  and \ \ f(n+1, m)=f(n, m+m^{2}).$$

By [4, Lemma  12.2],  $f(n,m)\leq 2^{2^{n-1}-1}   m^{2^{n-1}}.$

We use induction on $n$ to prove the following auxiliary statement
\vspace{3 mm}

$  \ \   (*)  $   {\it   Let  $ F, A_{1}, A_{2}, …, A_{n}$ be subsets of a semigroup  $S$ such that  $S=F^{-1}(A_{1}\bigcup A_{2}\bigcup …\bigcup A_{n})$
and $|F |\leq m$.   Then there exist $i\in  \{1,2,…,n\}$  and a subset $K$  of $S$ such that $S=K^{-1} \Delta (A_{i})$
and $|K|\leq  f(n,m).$}
\vspace{3 mm}

For $n=1,$ we have $S=F^{-1}  A_{1}$.  We take an arbitrary $x\in  S$  and choose $g\in  F$ such that   $xA_{1} \bigcap g^{-1} A_{1} \neq\emptyset.$  Then $A_{1} \bigcap x^{-1}g^{-1} A_{1} \neq\emptyset$, $A_{1} \bigcap (gx)^{-1}  A_{1}\neq\emptyset$ so $gx \in \Delta(A)$,  $x\in g^{-1} \Delta(A)$,  $x\in F^{-1} \Delta (A)$  and $S=F^{-1} \Delta (A)$.

Let  $S=F^{-1}(A_{1}\bigcup  A_{2} \bigcup  …  \bigcup A _{n+1}).$  We consider two cases.
\vspace{3 mm}

{\it Case 1. }  $gA_{1} \subseteq  F^{-1} (A_{2}\bigcup  … \bigcup A  _{n+1})$  for some $g \in S$.  Then   $A_{1} \subseteq  g^{-1} F^{-1} (A_{2}  … A _{n+1})$  and $ S=(F^{-1} \bigcup F^{-1} g^{-1}  F^{-1})(A_{2} \bigcup …\bigcup A _{n+1}).$  Since $\mid F\bigcup  FgF \mid\leq   m + m^{2}$,
 by the inductive hypothesis, there exist $i\in  \{2,3,…, n+1\}$  and a subset $K$  of $S$  such that
 $$ S= K^{-1}  \Delta (A_{i}),  \ \  |K|\leq  f (n,  m+m^{2})= f(n+1, m).$$

{\it Case 2. } $xA_{1} \bigcap  F^{-1} A_{1}  \neq\emptyset$  for every  $x\in S.$   Then $A_{1} \bigcap x^{-1} F^{-1} A_{1}\neq\emptyset$,  $x^{-1} g^{-1} A_{1} \bigcap A_{1}\neq\emptyset$  for some $g \in F$,  $gx \in\Delta  (A_{1})$  and  $x \in g^{-1} \Delta (A)$,  $S= F^{-1} \Delta  (A)$.  We set $K=F$  and note  that $|K |\leq m \leq f(n+1, m)$.

To conclude the proof, we assume that $S=A_{1}\bigcup ...\bigcup  A_{n}$, take an arbitrary $g\in S$, put $F= \{ g\}$, note that $S=F^{-1}(A_{1} \bigcup ...\bigcup A_{n})$  and apply $(*)$.
\vspace{3 mm}

{\it Proof of Theorem 2.}   Let $S$ be a finite  semigroup and $S= A_{1} \bigcup … \bigcup A_{n} .$  We take a minimal right ideal $R$ of $S$,  choose $r\in  S$ and note that $rS\subseteq  R$,   $S \subseteq r^{-1} R$,  so we may suppose that $S=R$. By [2, Theorem 1.63(g)],  $S$ is a direct product of a group $G$  and a right zero semigroup $I$.  We take $a \in I$  and put $H=G\times  \{a\}$.  For each $i \in \{1, …, n\}$, we denote $B_{i}=A_{i} \bigcap H$.  Since $H$  is a finite group, there are $j\in  \{1, …, n\}$ and $K\in H$  such that $|K|\leq n$  and  $H=K^{-1} \Delta_{H}(B_{j})$, where $\Delta_{H} (B_{j})=\{  x\in H: xB_{j} \bigcap B_{j}\neq\emptyset  \}  $.  We take an arbitrary $(g,b)\in  G \times I$  and choose $z \in K$ such that $z(g,a)\in \{x \in H: xB_{j} \bigcap B_{j}  \neq\emptyset\} $. Since $I$ is a right zero semigroup, we have $z(g,b)  B_{j} \bigcap B_{j}\neq\emptyset$.  Hence  $(g,b) \in  z^{-1} \{ x \in S: x A_{j} \bigcap  A_{j}\neq\emptyset\}$  and $ S=K^{-1 } \Delta (A_{j})$.
\vspace{3 mm}

{\it Proof of Theorem  3. }If $a \in A$  is left zero then, for every $x\in S$, we have  $S= a^{-1} a = a^{-1} \{ x \in A: xA\bigcap A \neq\emptyset\}$  and $S=a^{-1} \Delta(A)$.

If $a\in A$  is right zero then, for every $x \in S$, $a \in xA \bigcap A$ so  $xA \bigcap A\neq
\emptyset$  and $ S= \Delta (A)$  and $S= g^{-1} \Delta  (A)$ for each  $g \in S$.
\vspace{3 mm}

{\bf  Acknowledgement.} We thank George Bergman for constructive remarks on the seminal version of this note.
 \vspace{3 mm}

 \centerline{\bf References }
\vspace{3 mm}

[1]  Banakh, T., Protasov, I., Slobodianiuk, S. {\it  Densities, submeasures and partitions of $G$-spaces and groups},  Algebra and Discrete Mathematics, {\bf 17} (2014), Number 2, 193-221, preprint (http://arxiv.org./abs/1303.4612).

[2]  N. Hindman, D. Strauss, {\it  Algebra in the Stone-$\check{C}$ech Compactification}, 2nd edition, de Gruyter, 2012.

[3]  V.D. Mazurov, E.I. Khukhro (eds), {\it  Unsolved problems in group theory, the Kourovka notebook}, 13-th augmented edition, Novosibirsk, 1995.

[4]  Protasov I., Banakh T.,  {\it Ball structures and colorings of groups and graphs}, Math. Stud. Monogr. Ser., Vol.{\bf  11}, VNTL, Lviv, 2003.

\vspace{3 mm}
Department of Cybernetics, Kyiv University,

Prospect Glushkova 2, corp. 6,

03680 Kyiv, Ukraine

e-mail: I.V.Protasov@gmail.com

\vspace{3 mm}
Ksenia Protasova

Department of Cybernetics, Kyiv University,

Prospect Glushkova 2, corp. 6,

03680 Kyiv, Ukraine

e-mail: ksuha@freenet.com.ua

\end{document}